\documentclass[12pt]{article}
\usepackage{amsmath}
\usepackage{amssymb}
\newtheorem{thm}{Theorem}
\newtheorem{lemma}[thm]{Lemma}

\newtheorem{definition}[thm]{Definition}

\newcommand{\R}{\mathbb{R}}
\newcommand{\Rn}{\mathbb{R}^n}
\newcommand{\dif}[0]{\ensuremath{\,\mathrm{d}}}

\newcommand{\abs}[1]{\ensuremath{\vert #1 \vert}}

\DeclareMathOperator*{\diam}{diam}

\begin{document}
\title{STABILITY FOR THE INFINITY-LAPLACE EQUATION WITH VARIABLE
  EXPONENT}
\author{Erik Lindgren \qquad Peter Lindqvist}
\date{Norwegian University of Science and Technology}
\maketitle

{\small \textsc{Abstract:} \textsf{We study the stability for the
    viscosity solutions of the differential equation}
$$ \sum u_{x_i}u_{x_j}u_{x_i x_j}+\abs{\nabla u}^2\ln(\abs{\nabla
      u})\langle\nabla u, \nabla \ln p \rangle=0$$ \textsf{under
      perturbations of the function} $p(x).$ \textsf{The differential operator
  is the so-called $\infty(x)$-Laplacian.}}

\section{Introduction}

The object of our study is the curious differential equation 
\begin{equation}
  \label{eq:curious}
  \sum_{i,j=1}^{n} u_{x_i}u_{x_j}u_{x_i x_j}+\abs{\nabla u}^2\ln(\abs{\nabla
      u})\langle\nabla u, \nabla \ln p \rangle=0
\end{equation}
in a bounded domain $\Omega$ in $\R^n$. Here $p(x)$ is a positive
function, the so-called variable exponent, and it is of class
$C^1(\overline{\Omega})$. The equation comes from the mini-max problem
of determing
\begin{displaymath}
  \min_u\max_x\left\{\abs{\nabla u(x)}^{p(x)}\right\}.
\end{displaymath}
The case of a constant  $p(x)=p$ reduces to the celebrated
$\infty$-Laplace equation
\begin{equation}
  \label{eq:infty-lap}
  \Delta_\infty u\equiv\sum_{i,j=1}^{n} u_{x_i} u_{x_j}u_{x_i x_j}=0
\end{equation}
found by G. Aronsson. In [LL] the equation was derived as the limit
of the Euler-Lagrange equations for the variational integrals 
\begin{equation*}\label{eq:kpx-integrals}
  \left\{\int_\Omega\abs{\nabla u(x)}^{kp(x)}
    \dif x \right\}^{\frac{1}{k}}
\end{equation*}
as $k\to\infty$. Such integrals were first considered by Zhikov,
cf. [Z]. See also [RMU] for similar equations. For sufficiently smooth
solutions the meaning of the equation is that 
$$\abs{\nabla u(x)}^{p(x)} = C$$
along any fixed stream line (different stream lines may have different
constants attached). ---In general, solutions have to be interpreted
in the \emph{viscosity sense}, and we assume that the reader is
acquainted with the basic theory of viscosity solutions, see [CIL, K,
C].

The viscosity solution with prescribed Lipschitz continuous boundary
values is unique, cf. [LL]. Taking into account that, in contrast,
uniqueness does not always hold for the $\infty$-Poisson equation
$$\Delta_{\infty}u = \varepsilon(x),$$
as an example with a uniformly continuous sign-changing function $\varepsilon(x)$ in [LW]
shows, the uniqueness for the curious equation \eqref{eq:curious} is
pretty remarkable. Therefore we have found it worth our while to study
the stability under variations of $p(x).$

Our first result is about a perturbation of the $\infty$-Laplace
equation
\eqref{eq:infty-lap}.

\begin{thm}
Let $p \in C^{1}(\overline{\Omega})$ be a positive function and
suppose that $u \in C(\overline{\Omega})$ is the viscosity solution of
$$\Delta_{\infty}u + |\nabla u|^{2}\ln(|\nabla u|)\langle \nabla u,
\nabla \ln(p) \rangle = 0$$
and that $v \in C(\overline{\Omega})$ is the viscosity solution of
$$\Delta_{\infty}v = 0,$$
both having the same Lipschitz continuous boundary values $f$. Then
the estimate
\begin{equation}
\|v-u\|_{L^{\infty}(\Omega)} \leq C_{1} \|\nabla \ln p
\|_{L^{\infty}(\Omega)} + C_{2} \|\nabla \ln p
\|_{L^{\infty}(\Omega)}^{\frac{1}{5}}\vert\ln(c \|\nabla \ln p
\|_{L^{\infty}(\Omega)})\vert
\end{equation}
is valid with constants depending only on
$\|f\|_{W^{1,\infty}(\Omega)}$ and $\diam(\Omega).$
\end{thm}

An interpretation is that when $p(x)$ deviates only little from a
constant value, then $u$ is close to $v$. But, as we have pointed out,
the perturbation term $ |\nabla u|^{2}\ln(|\nabla u|)\langle \nabla u,
\nabla \ln(p) \rangle $ cannot be replaced by an arbitrary small
perturbation $\varepsilon(x)$, despite the possibility to select
$p(x)$ in any manner. ---The exponent $\frac{1}{5}$ seems to be an
artifact of the arrangements in our proof in section 4.

We also address the problem with two  positive exponents $p_1, p_2 \in
C^1(\overline{\Omega}),$ but now the result is
weaker. Suppose that $u_{\nu} \in C(\overline{\Omega})$ is a viscosity
solution of
\begin{equation}
\Delta_{\infty}u_{\nu} + |\nabla u_{\nu}|^{2}\ln(|\nabla u_{\nu}|)\langle \nabla u_{\nu},
\nabla \ln(p_{\nu}) \rangle = 0
\end{equation}
in $\Omega$, $\nu = 1,2.$ If $u_1$ and $u_2$ have the same Lipschitz
continuous boundary values, then
\begin{gather}
\|u_1 - u_2 \|_{L^{\infty}(\Omega)} \,\leq \, \frac{\text{Const.}}{\vert
  \ln \bigl(\|\nabla \ln p_2 -
  \nabla \ln p_1\|_{L^{\infty}(\Omega)}\bigr)|^{\kappa}},
\end{gather}
where $\kappa > 0$ depends on $\max(p_{\nu}), \min(p_{\nu}).$ The
constant depends on the boundary values and on the norms $\|\nabla \ln
p_{\nu}\|_{\infty}.$ 
Needless to say, the obtained modulus of stability  appears to
be far from sharp. Therefore we have only sketched out the proof in section
5. In the one-dimensional case a sharp bound is easily reached via the
``first integral''
$|u'_{\nu}(x)|^{p_{\nu}(x)} = C_{\nu}.$

\section{Preliminaries}

We briefly recall some basic concepts. Let $\Omega$ be a bounded
domain in $\Rn$ and suppose that $f:\,\partial \Omega \rightarrow \R$
is a Lipschitz continuous function satisfying
$$|f(x) - f(y)| \leq L |x - y|.$$
By extension, we may as well assume that the inequality holds in the
whole space, if needed. The abbreviation 
$$\Delta_{\infty(x)} u \equiv \Delta_{\infty} u + |\nabla u|^2
\ln(|\nabla u|)\langle \nabla u, \nabla \ln p \rangle $$
is convenient\footnote{The suggestive subscript $\infty(x)$ symbolizes
  the ``variable exponent infinity''.}.

 To be on the safe side, we assume that $p\in
C^1(\overline{\Omega})$, $p(x)>0$. Then viscosity solutions
to the equation \eqref{eq:curious} can be defined in the standard way.

\begin{definition}\label{def:viscosity}
  We say that a lower semicontinuous function $v:\Omega\to
  (-\infty,\infty]$ is a \emph{viscosity supersolution} if, whenever
  $x_0\in \Omega$ and $\varphi\in C^2(\Omega)$ are such that
  \begin{enumerate}
  \item $\varphi(x_0)=u(x_0)$, and
  \item $\varphi(x)<v(x)$, when $x\not=x_0$, 
  \end{enumerate}
  then we have
  \begin{displaymath}
    \Delta_{\infty(x_0)}\varphi(x_0)\leq 0.
  \end{displaymath}
\end{definition}

The \emph{viscosity subsolutions} have a similar definition; they are
upper semicontinuous, the test functions touch from above and the
differential inequality is reversed. Finally, a \emph{viscosity
  solution}  is both a viscosity supersolution and
viscosity subsolution. 

There is an alternative way of expressing the definition in terms of
''semijets''. 

\section{Auxiliary Equations}

Following R. Jensen in [J]  we introduce two auxiliary equations. For a
constant exponent $p$ the situation is
\begin{align}
  &\max\{\varepsilon-\abs{\nabla u^{+}},\Delta_{\infty} u^{+}\}=0 
  & &\text{\textit{Upper equation}}\label{eq:upper}\\
  &\Delta_{\infty}u=0 &
  &\text{\textit{Equation}} \nonumber
\\ 
  &\min\{\abs{\nabla u^{-}}-\varepsilon,\Delta_{\infty}u^{-}\}=0   
  & &\text{\textit{Lower equation}}\label{eq:lower}
\end{align}
where $\varepsilon > 0.$ Given  $\varepsilon > 0$, three viscosity
solutions
$u^{-}, u, u^{+}$ are constructed with  the same boundary values $f$
so that
\begin{gather}
u^{-} \leq u\leq  u^{+}\nonumber\\
\| u^{+} -u^{-}\|_{L^{\infty}(\Omega)} \leq \varepsilon \diam(\Omega)\\
\|\nabla u^{\pm}\|_{L^{\infty}(\Omega)} \leq K + \varepsilon =
 K_{\varepsilon}\nonumber
\end{gather}
where $K$ depends only on the Lipschitz constant $L$ of $f$. The
virtue of the auxiliary equations is that 
$$ \varepsilon-  |\nabla u^{+}| \leq 0,\quad  |\nabla u^{-}|-
\varepsilon \geq 0$$
in the viscosity sense. We refer to [J] and [LL] about the
construction via variational integrals.

We need a \emph{strict} supersolution. We will construct a function
$g(u^{+}) \approx u^{+}$ such that $\Delta_{\infty}g(u^{+}) < 0.$ To this end
we use the following \emph{approximation of the identity}
\begin{equation}
g(t) = \frac{1}{\alpha}\ln \bigl(1+A(e^{\alpha t}-1)\bigr), \quad A >
1,  \alpha > 0
\end{equation}
taken from [JLM] and [LL]. For $t>0, A > 1, \alpha > 0$ we have
\begin{gather*}
0 < g(t) -1 < \frac{A-1}{\alpha}\\
(A-1) e^{-\alpha t} < g'(t) < A-1\\
\frac{g''(t)}{g'(t)} = - \alpha\bigl(g'(t)-1\bigr),
\end{gather*}
which are easy to verify.

\begin{lemma} Let $v > 0$ and consider $w = g(v).$ If
$$ \varepsilon - |\nabla v| \leq 0 \quad\text{and}\quad \Delta_{\infty}v
\leq 0$$
in the viscosity sense, then the inequality
\begin{equation}
\Delta_{\infty} w \leq - \alpha
(A-1)A^{-1}e^{-\alpha\|v\|_{\infty}}\varepsilon^{4} \equiv -\mu
\end{equation}
holds in the viscosity sense.
\end{lemma}

\emph{Proof:} Formally, the equation for  $w = g(v)$ is
\begin{align*}
\Delta_{\infty} w &= g'(v)^3\Delta_{\infty}v + g''(v)g'(v)^{2}|\nabla
v|^{4}\\
& \leq 0 +  g''(v)g'(v)^{2}|\nabla
v|^{4}\\
&= - \alpha (A-1)A^{-1}e^{-\alpha v}g'(v)^4|\nabla v|^4\\
&\leq  - \alpha (A-1)A^{-1}e^{-\alpha v}1^4\varepsilon^4.
\end{align*}
To conclude the proof, one has to pass the calculation over to test
functions. $\Box$

\medskip
We will apply the lemma on $w = g(u^+)$ and we assume that $f > 0$ so
that the encountered functions are non-negative. It holds that
\begin{equation}
\min_{\partial \Omega}(f) = \min_{\partial \Omega}(u) \leq u \leq u^+
\leq \max_{\partial \Omega}(f) + \varepsilon \diam(\Omega)
\end{equation}
by the maximum principle and (8). Fix
$$ \alpha = \frac{1}{\|u^+\|_{\infty}}.$$
Estimate (10) in the lemma above becomes
\begin{equation}
\label{eq:strict}
\Delta_{\infty}g(u^+) \leq - \mu = -
\frac{(A-1)\varepsilon^4}{Ae\|u^+\|_{\infty}}.
\end{equation}

\section{Proof of the Stability}

Suppose that $u_{1}$ is a viscosity (sub)solution of
$$\Delta_{\infty}u_1 + |\nabla u_1|^{2}\ln(|\nabla u_1|)\langle \nabla
u_1,\nabla \ln p_1 \rangle = 0$$ 
and that $u_2$ is a viscosity (super)solution of
$$\Delta_{\infty}u_2 = 0, \qquad (p_2 = \text{constant})$$
both with boundary values $f$. Adding the same constant to $f, u_1,$
and $u_2$, we may assume that $f \geq 0$ and $u_2 \geq 0.$ Given
$\varepsilon > 0,$ write
$$v_2 = u_2^{+},\quad w_2 = g(v_2) = g(u_2^+).$$
We obtain the estimate
\begin{align}
u_1-u_2 &= (u_1-w_2) + (w_2-v_2) + (v_2-u_2)\nonumber\\
&<  (u_1-w_2) +\frac{A-1}{\alpha} + \varepsilon \diam(\Omega).
\end{align}
The last two terms could be made as small as we please, but the term
$u_1-w_2$ requires our attention, since there $w_2$ depends also on $A$ and
$\varepsilon.$ 

\begin{lemma}
We have
\begin{equation}
u_1-w_2 \leq C_{\varepsilon}^3\,\|u_2^+\|^2_{\infty}\,\varepsilon^{-4}\ln\Bigl(
\frac{C_{\varepsilon}}{\varepsilon}\Bigr)\,\|\nabla \ln p_1\|_{\infty},
\end{equation}
where $ C_{\varepsilon} = C(1+\varepsilon).$
\end{lemma}

\emph{Proof:} Let $\sigma = \max(u_1-w_2).$ If $\sigma \leq 0$ there
is nothing to prove. Assume thus that $\sigma > 0.$ In order to use
the Theorem on Sums for viscosity solutions, we double the variables
writing
\begin{equation}
\label{eq:M}
    M_j=\sup_{\substack{x\in\Omega\\y\in\Omega}}\left(u_1(x)-w_2(y)
      -\frac{j}{2}\abs{x-y}^2\right)
  \end{equation}
as usual. Then $M_j \geq \sigma$ (take $x = y$ to see this). 
  The supremum is attained at some points $x_j$, $y_j$. Now $|x_j-y_j|
  \to 0$ as $j \to \infty$ and 
  \begin{displaymath}
    x_j\to \Hat{x}, \qquad y_j\to \Hat{y} = \Hat{x},
  \end{displaymath}
at least for a subsequence. We claim that
 $\Hat{x}$ is an interior point of $\Omega.$ Indeed, if $\Hat{x}
 \in \partial \Omega$ then
\begin{align*}
u_1(\Hat{x})- w_2(\Hat{x})&=(u_1(\Hat{x})- v_2(\Hat{x})) + (
v_2(\Hat{x})- w_2(\Hat{x}))\\
&= 0 +  (v_2(\Hat{x})-g(v_2(\Hat{x})) \leq 0
\end{align*}
and hence 
\begin{align*}
u_1(\Hat{x})- w_2(y) &= (u_1(\Hat{x})- w_2(\Hat{x}))+
(w_2(\Hat{x})-w_2(y))\\
&\leq w_2(\Hat{x})-w_2(y),
\end{align*}
which, by continuity, is less than $\sigma/2 < \sigma \leq M_j$
provided that $|\Hat{x}-y|$ is small. Hence $\Hat{x} \in \Omega.$

We conclude that also $x_j$ and $y_j$ are interior points for large
indices $j.$ We need the bounds
\begin{equation}
\varepsilon \leq j|x_j-y_j| \leq C_{\varepsilon}.
\end{equation}
The upper bound follows from
\begin{gather*}
 u_1(x_j)-w_2(y_j)-\frac{j}{2}\abs{x_j-y_j}^2 \geq
u_1(x_j)-w_2(x_j),\\
\frac{j}{2}\abs{x_j-y_j}^2 \leq w_2(x_j)-w_2(y_j)
\leq \|g'(v_2)\nabla v_2\|_{\infty}|x_j-y_j|
\leq AK_{\varepsilon}|x_j-y_j|,
\end{gather*}
where we used that $g'(v_2) < A.$ We had $K_{\varepsilon} = K +
\varepsilon$ and we will later see that $A \leq 2.$ Then
$C_{\varepsilon} = 2 K_{\varepsilon}$ will do. The lower bound is
deduced from the fact that $\varepsilon - |\nabla w_2| \leq 0$ in the
viscosity sense ($\nabla w_2 = g'(v_2)\nabla v_2,\, 1 \leq g'(v_2),\,
\varepsilon \leq |\nabla v_2|$). To wit,
  \begin{displaymath}
    u_1(x_j)-w_2(y_j)-\frac{j}{2}\abs{x_j-y_j}^2 \geq 
    u_1(x_j)-w_2(y)-\frac{j}{2}\abs{x_j-y}^2,
  \end{displaymath}
from which it follows that the function
$$\varphi(y) =
w_2(y_j)+\frac{j}{2}\abs{x_j-y_j}^2-\frac{j}{2}\abs{x_j-y}^2$$
touches $w_2(y)$ from below at the point $y_j.$ Thus $\varepsilon \leq
|\nabla \varphi(y_j)|,$ and this is the desired inequality, indeed.

According to the Theorem on Sums there exist symmetric $n\times
n$-matrices  $\mathbb{X}_j$ and $\mathbb{Y}_j$ such that $\mathbb{X}_j\leq \mathbb{Y}_j$ and
  \begin{align*}
    (j(x_j-y_j),\mathbb{X}_j)\in & \overline{J^{2,+}}u_1(x_j),\\
    (j(x_j-y_j),\mathbb{Y}_j)\in & \overline{J^{2,-}}w_2(y_j)
  \end{align*}
  where $\overline{J^{2,+}}u_1(x_j)$ and $\overline{J^{2,-}}w_2(y_j)$ are
  the closures of the super- and subjets. (Caution:
  \emph{super}solutions are tested with \emph{sub}jets.) For the jets
  and their closures we refer to [CIL], [C], [K]. The meaning
  of the notion is that we can rewrite
  the equations as
  \begin{align*}
    j^2&\langle \mathbb{Y}_j(x_j-y_j),x_j-y_j\rangle  \leq -\mu,\\
     j^2&\langle \mathbb{X}_j(x_j-y_j),x_j-y_j\rangle \\
    &+ j^3\abs{x_j-y_j}^2 \ln(j\abs{x_j-y_j})\langle x_j-y_j,
    \nabla \ln p_1(x_j)\rangle \geq 0,\\
    j&\abs{x_j-y_j}\geq \varepsilon,\\
    j&\abs{x_j-y_j}\leq C_{\varepsilon}.
  \end{align*}
It follows that
\begin{gather*}
 j^2\langle(\mathbb{Y}_j-\mathbb{X}_j)(x_j-y_j),x_j-y_j\rangle  \leq -\mu
+ C_{\varepsilon}^{3}\ln \bigl(\frac{
  C_{\varepsilon}}{\varepsilon}\bigr)\|\nabla \ln p_1\|_{\infty}.
\end{gather*}
The left-hand member is a positive semidefinite quadratic form, since \\
$\mathbb{Y}_j-\mathbb{X}_j \geq 0,$ in other words it is non-negative. Thus
$$ \mu \leq C_{\varepsilon}^{3}\ln \bigl(\frac{
  C_{\varepsilon}}{\varepsilon}\bigr)\|\nabla \ln p_1\|_{\infty}.$$
Recall the expression for $\mu$ in \eqref{eq:strict}. The above estimate can be
written as
$$\frac{(A-1)\varepsilon^{4}}{Ae\,\|u_2^{+}\|_{\infty}}
\leq C_{\varepsilon}^{3}\ln \Bigl(\frac{
  C_{\varepsilon}}{\varepsilon}\Bigr)\|\nabla \ln p_1\|_{\infty}.$$
We fix $A > 1$ so that
$$\frac{A-1}{\alpha} = \sigma,$$
where we had $\alpha^{-1} = \|u_2^+\|_{\infty}.$ Then
\begin{equation*}
\sigma \leq Ae\varepsilon^{-4} C_{\varepsilon}^{3}\ln \bigl(\frac{
  C_{\varepsilon}}{\varepsilon}\bigr)\|\nabla \ln
p_1\|_{\infty}\,\|u_2^{+}\|^{2}_{\infty}.
\end{equation*}
Further, $A \leq 2$ so that $Ae$ can be absorbed into the constant $
C_{\varepsilon}.$ (Indeed, $\sigma = \max(u_1-w_2) \leq \max(u_1) =
u_1(\xi)$ where $\xi$ is some boundary point. Now
$$A = 1+\alpha \sigma \leq 1 +
\frac{u_1(\xi)}{\|u_2^{+}\|_{\infty}} \leq  1 +
\frac{u_1(\xi)}{u_2^{+}(\xi)} = 2,$$
because the functions have the same boundary values.) This concludes
the proof of the estimate (14). $\Box$

\medskip
We return to (13). Using (14) we obtain
\begin{gather*}
u_1-u_2 \leq \sigma + \sigma + \varepsilon \diam(\Omega)\\
\leq 2 \varepsilon^{-4} C_{\varepsilon}^{3}\ln \bigl(\frac{
  C_{\varepsilon}}{\varepsilon}\bigr)\|u_2^{+}\|^{2}_{\infty} \,\|\nabla \ln
p_1\|_{\infty} + \varepsilon \diam(\Omega).
\end{gather*}
It remains to determine $\varepsilon$ nearly optimally. To simplify, we use
(11):
$$u_2^+ \leq u_2 +  \varepsilon \diam(\Omega) \leq \|f\|_{\infty}  +
\varepsilon \diam(\Omega).$$
We have to optimize 
$$(C+\varepsilon)^{3}\ln\bigl(\frac{C+\varepsilon}{\varepsilon}\bigr)
(\|f\|_{\infty}+ \varepsilon)^{2}\|\nabla \ln p_1\|_{\infty}\,
\varepsilon^{-4} +  \varepsilon \diam(\Omega),$$
which, renaming constants, is the same as an expression of the form
$$\Bigl(\frac{C+\varepsilon}{\varepsilon}\Bigr)^{5}
\ln\Bigl(\frac{C+\varepsilon}{\varepsilon}\Bigr)\|\nabla \ln
p_1\|_{\infty}\,\varepsilon + \varepsilon a.$$
We consider two cases. The case of a large $\|\nabla \ln
p_1\|_{\infty}$
 is plain.
 Namely, if $a \leq 32 \,\|\nabla \ln p_1\|_{\infty}$ we just take
$\varepsilon = 1$ and obtain immediately a majorant of the form
 $C_1 \|\nabla \ln p_1\|_{\infty}$. If not, we can determine
 $\varepsilon$ from the equation
$$\Bigl(\frac{C+\varepsilon}{\varepsilon}\Bigr)^{5}\|\nabla \ln
p_1\|_{\infty} = a.$$
This yields a majorant like
$$C_2\|\nabla \ln p_1\|_{\infty}^{\frac{1}{5}}\ln(c\,\|\nabla \ln
p_1\|_{\infty}).$$
Combining the two cases we arrive at the desired estimate (3), yet so far
only for $\max(u_1-u_2).$ The corresponding estimate for $\max(u_2-u_1)$ is still missing;
the situation is not symmetric. To complete the proof, observe that
$$u_2(x)-u_1(x) = (k-u_1(x))-(k-u_2(x))$$
where the constant is large enough to make the new viscosity
solution $k-u_2$ positive; $k = \max(f)$ will do. Now
$\Delta_{\infty}(k-u_2(x)) = 0$ and the situation has been reduced to
the previous case. ---Instead, we could have repeated the proof, this
time using the Lower Equation (7).  $\Box$

\section{Two Varying Exponents}

In the case of two exponents none of which is constant, an extra
complication arises: the parameter $\alpha$ must be taken very large,
say $\alpha \approx \varepsilon^{-1},$ and then the exponential factor
in the counterpart to \eqref{eq:strict} is extremely small. This
weakens the final result.

In principle, the proof is a repetition of the previous one. Only an
outline is provided below. First, the auxiliary equations in section 3 are
modified so that $\Delta_{\infty}$ is replaced by
$\Delta_{\infty(x)}.$ As in [LL] one then obtains the estimate
\begin{equation}
 \label{eq:kappa}
\|u^{+}-u^{-}\|_{\infty} \leq B \varepsilon^{\kappa}
\end{equation}
where $\kappa > 0$ (either $\kappa = \min(p(x))$ or $\kappa =
\max(p(x))$).
Second, we need a strict supersolution to equation (1).

\begin{lemma} Consider $w = g(v)$ for $v > 0.$ If
$$\varepsilon - |\nabla v| \leq 0\quad \text{and}\quad
\Delta_{\infty(x)}v \leq 0$$
in the viscosity sense, then
\begin{equation}
\label{eq:strictx}
\Delta_{\infty(x)} w \leq - \varepsilon^{3}(A-1)A^{-1}e^{-(1+\|\nabla
  \ln p\|_{\infty})\|v\|_{\infty}\varepsilon^{-1}}\equiv -\mu
\end{equation}
in the viscosity sense.
\end{lemma}

\emph{Proof:} A routine calculation yields
\begin{align*}
\Delta_{\infty(x)} w &\leq g'(v)^3(g'(v)-1)|\nabla
v|^3\left\{-\alpha|\nabla v|+|\nabla\ln p|\right\}\\
 &\leq g'(v)^3(g'(v)-1)|\nabla
v|^3\left\{-\alpha \varepsilon+\|\nabla\ln p\|_{\infty}\right\}.
\end{align*}
Given $\varepsilon > 0$, we fix $\alpha = \alpha(\varepsilon)$ so that
$$-\alpha\varepsilon + \|\nabla\ln p\|_{\infty} = -1.$$
The estimate \eqref{eq:strictx} readily follows. $\Box$

\medskip
Suppose now that $u_{\nu}$ is a viscosity solution of the equation
(4), $\nu = 1,2.$ We assume that $u_1 = u_2 = f$ on $\partial \Omega.$
By adding a constant, we reach the situation that $u_2^+ \geq u_2 >
0.$ Write
$$v_2 = u_2^+,\quad w_2 = g(v_2).$$
Now
\begin{align*}
u_1-u_2 &=(u_1-w_2) + (w_2-v_2) + (v_2-u_2)\\
        &\leq(u_1-w_2)  + \frac{A-1}{\alpha} + B\varepsilon^{\kappa}.
\end{align*}

\begin{lemma} We have
\begin{equation}
u_1-w_2
\leq 2 \varepsilon^{-2}A C_{\varepsilon}^{3}\ln
\Bigl(\frac{C_{\varepsilon}}{\varepsilon}\Bigr) e^{(1+\|\nabla
  \ln p_2\|_{\infty})\,\|v_2\|_{\infty}\,\varepsilon^{-1}}\|\nabla \ln p_2
- \nabla \ln p_1\|_{\infty}.
\end{equation}
\end{lemma}

\emph{Proof:} Denote $\sigma =  \max(u_1-u_2).$ We may assume that
$\sigma > 0.$ Double the variables as in \eqref{eq:M}. By the Theorem
on Sums we again obtain symmetric matrices $\mathbb{X}_j$ and
$\mathbb{Y}_j$ so that $\mathbb{X}_j\leq \mathbb{Y}_j$ and
  \begin{align*}    
    j^2&\langle \mathbb{Y}_j(x_j-y_j),x_j-y_j\rangle \\
    &+ j^3\abs{x_j-y_j}^2 \ln(j\abs{x_j-y_j})\langle x_j-y_j,
    \nabla \ln p_2(y_j)\rangle \leq -\mu,\\
     j^2&\langle \mathbb{X}_j(x_j-y_j),x_j-y_j\rangle \\
    &+ j^3\abs{x_j-y_j}^2 \ln(j\abs{x_j-y_j})\langle x_j-y_j,
    \nabla \ln p_1(x_j)\rangle \geq 0,\\ 
    j&\abs{x_j-y_j}\geq \varepsilon,\\
    j&\abs{x_j-y_j}\leq C_{\varepsilon}.
  \end{align*}
Write $\ln p_1(x_j) = \ln p_1(x_j) - \ln p_2(x_j) + \ln p_2(x_j)$ and
arrange the equations. It follows that
  \begin{align*}
    0 &\leq j^2\langle (\mathbb{Y}_j-\mathbb{X}_j)(x_j-y_j),(x_j-y_j)\rangle 
         \leq   -\mu\\ &+j^3\abs{x_j-y_j}^2\ln(j\abs{x_j-y_j})
      \langle x_j-y_j,\nabla \ln
      p_2(x_j)-\nabla\ln p_2(y_j)\rangle\\
&+j^3\abs{x_j-y_j}^2\ln(j\abs{x_j-y_j})
      \langle x_j-y_j,\nabla \ln
      p_1(x_j)-\nabla\ln p_2(x_j)\rangle\\
     & \leq -\mu+ C^3_{\varepsilon}\ln\Bigl(\frac{C_{\varepsilon}}{\varepsilon}\Bigr)
      \abs{\nabla \ln p_2(x_j)-\nabla\ln p_2(y_j)}\\
&+ C^3_{\varepsilon}\ln\Bigl(\frac{C_{\varepsilon}}{\varepsilon}\Bigr)
      \|\nabla \ln p_1- \nabla\ln p_2\|_{\infty}.
    \end{align*}
As $j \to \infty$, $x_j-y_j \to 0.$ so that by continuity
\begin{equation}
\mu \leq  C^3_{\varepsilon}\ln\Bigl(\frac{C_{\varepsilon}}{\varepsilon}\Bigr)
      \|\nabla \ln p_1- \nabla\ln p_2\|_{\infty}.
\end{equation}
To conclude the proof, we fix $A > 1$ so that
$$\frac{A-1}{\alpha} = \frac{\sigma}{2}$$
and insert the expression for $\mu$ given in \eqref{eq:strictx}. Hence 
\begin{gather*}
\varepsilon^{3}\frac{\sigma \alpha}{2}A^{-1}e^{-(1+\|\nabla
  \ln p_2\|_{\infty})\|v_2\|_{\infty}\varepsilon^{-1}}\\
\leq  C^3_{\varepsilon}\ln\Bigl(\frac{C_{\varepsilon}}{\varepsilon}\Bigr)
      \|\nabla \ln p_1- \nabla\ln p_2\|_{\infty}.
\end{gather*}
The estimate follows, since $\alpha \approx 1/\varepsilon.$ $\Box$

\medskip
In order to finish the proof of (4) we choose $\varepsilon$ in the inequality
\begin{multline*}
u_1-u_2 \\ \leq 4 A \varepsilon^{-2} C^3_{\varepsilon}\ln\Bigl(\frac{C_{\varepsilon}}{\varepsilon}\Bigr)e^{(1+\|\nabla
  \ln p_2\|_{\infty})\|v_2\|_{\infty}\,\varepsilon^{-1}}
\|\nabla \ln p_1- \nabla\ln p_2\|_{\infty}
+ B \varepsilon^{\kappa}
\end{multline*}
so that
$$e^{(1+\|\nabla
  \ln p_2\|_{\infty})\|v_2\|_{\infty}\varepsilon^{-1}}
\|\nabla \ln p_1- \nabla\ln p_2\|_{\infty} \approx
\varepsilon^{\kappa+2}.$$
We omit the calculation.

\end{document}